%% file: paper1.tex
\documentclass[12pt,leqno]{article}
\usepackage{amscd,amsmath,amsthm,amssymb,psfig}
\theoremstyle{plain}
\newtheorem{X}{X}[section] 
\newtheorem{Thm}[X]{Theorem}
\newtheorem{Pptn}[X]{Proposition}
\newtheorem{Lem}[X]{Lemma}  
\newtheorem{Cor}[X]{Corollary}

\theoremstyle{definition}
   
\newtheorem{Def}[X]{Definition}
\newtheorem{Rem}[X]{Remark}

\newtheorem{E}{}

\renewcommand{\Bbb}{\mathbb}  

\def\cal{\mathcal}
\renewcommand{\frak}{\mathfrak}

\def\primes{p_1,\dots,p_t}

\def\hset{H(P)}
\def\hsetpr{H(P')}
\def\zzt{\big(\Bbb Z/3\Bbb Z\big)^t}

\def\fp{\frak p}
\def\idfac#1#2{\fp_1^{{#1}_1}(\bar\fp_1)^{{#2}_1}\dots\fp_t^{{#1}_t}
(\bar\fp_t)^{{#2}_t}}
\def\idfacsh#1{\fp_1^{{#1}_1}\dots\fp_t^{{#1}_t}}
\def\image{\operatorname{im}}
\def\cok{\operatorname{cok}}
\def\zh{\hat{\Bbb Z}}
\def\E{E:y^2=x^3+a(x-b)^2}
\def\Ep{E':y^2=x^3-27a(x-4a-27b)^2}
\def\pio{\pi_0(\overline E)}
\def\gal{\operatorname{Gal}}

\begin{document}

\title{A Formula for The Selmer Group of a Rational
Three-Isogeny}

\author{Matt DeLong}

\maketitle

\begin{abstract} A formula is given for the dimension of the Selmer group
of the rational three-isogeny of elliptic curves of the form
$y^2=x^3+a(x-b)^2$.  The formula is in terms of the three-ranks of the
quadratic number fields $\Bbb Q(\sqrt{a})$ and $\Bbb Q(\sqrt{-3a})$ and
various aspects of the arithmetic of these number fields. In addition a
duality theorem is used to relate the dimension of the Selmer group of the
three-isogeny with the dimension of the Selmer group of its dual isogeny.
\end{abstract} 

\input{introp1} 
\input{toppap}

\bibliographystyle{plain}
\bibliography{paper1}

DEPARTMENT OF MATHEMATICS, TAYLOR UNIVERSITY

236 W. READE AVE., UPLAND, IN 46989

{\it E-mail address:}  {\tt mtdelong@tayloru.edu}

\end{document}

%% file: introp1.tex
\section{Introduction}\label{S:intro}


We study elliptic curves which admit a rational 3-isogeny.  Such an
elliptic curve, $E/\Bbb Q$, has a subgroup $T\subset E(\bar\Bbb Q)$ of
order 3 defined over $\Bbb Q$.  We may suppose that $E$ is given by the
equation $y^2=x^3+ax^2+cx+d$, and by a change of coordinates, we may
assume that $T$ is generated by the point $(0,\sqrt d)$.  Using the
explicit addition law found in~\cite{jS92}, we find that this point has
order 3 precisely when $c^2=4ad$.  If $c\neq 0$, then the equation of the
curve can be written as $y^2=x^3+a(x-b)^2$.  

In ~\cite{jT93} Top 
gave an inequality
relating the rank of the Mordell-Weil group over $\Bbb Q$ of an elliptic
curve of the form $y^2=x^3+a(x-b)^2$ to the three-ranks of the quadratic
number fields $\Bbb Q(\sqrt a)$ and $\Bbb Q(\sqrt{-3a})$, denoted $r_3(a)$
and $r_3(-3a)$ respectively. Here
we expand on the work of Top to give a formula for the dimension of the
Selmer group of the three-isogeny from the elliptic curve given by the
equation $y^2=x^3+a(x-b)^2$ to the isogenous curve given by
$y^2=x^3-27a(x-(4a+27b))^2$. The formula relates the dimension of this
Selmer group to $r_3(a)$ and $r_3(-3a)$ as well as the cardinalities of
certain sets of primes dividing $a$ and $4a+27b$.  In addition, we
calculate via a formula of Cassels~\cite{jC65} the difference between the
dimensions of the Selmer groups of the isogeny and its dual. This provides
a useful check on the main result.

If $c=0$, then the equation of the curve can be written as
$y^2=x^3+d$.  We extend the work found in~\cite{pS86} and study elliptic
curves of this latter form in~\cite{mD98}. 

%% file: toppap.tex
\section{The Work of Top}\label{S:Top}

Let $E$ be the elliptic curve defined over $\Bbb Q$ given by the equation
$\E$ with $a,b\in\Bbb Z$ such that $a,b,\text{ and }4a+27b$ are non-zero. 
Up to an isomorphism over $Q$, we may assume that $a$ is square-free, and
so we do so throughout.  This curve admits a rational 3-isogeny.  We
denote this isogeny and its dual by $\psi$ and $\psi'$ respectively.

In~\cite{jT93} Top gave an
injection of $S^\psi$ into a group related to the class group of $\Bbb
Q(\sqrt{-3a})$, 
and used this to bound the dimension of the Selmer group.  He
similarly bounded the dimension of $S^{\psi'}$, and so deduced the following
theorem concerning the rank of an elliptic curve of the form $\E$. 

\begin{Thm}\label{T:Top} (Top)
Denote by $E_{a,b}/\Bbb Q$ the elliptic curve which is given by the
equation $y^2=x^3+a(x-b)^2$ with $a,b\in\Bbb Z,\,a\ne 0,\,b\ne
0,\,4a+27b\ne 0$.  Assume that 
\begin{align}
(1)&\, a\equiv 2 \text{ or } 11\pmod {12}\text{ and}\notag\\
(2)&\, a\text{ is square-free.}\notag 
\end{align}

Write s for the number of primes $p\ge 5$ such that $p|b$ and
$(\frac ap)=1$.  Similarly, 
write t for the number of primes $p\ge 5$ such that $p|(4a+27b)$ and
$(\frac{-3a}p)=1$.

Then
$$\operatorname{rank}E_{a,b}(\Bbb Q)\le r_3(a)+r_3(-3a)+s+t+1.$$
\end{Thm}

\begin{Rem} The first condition of the theorem is really two separate
conditions.  Condition (1a), $a\equiv 2 \text{ or } 3\pmod 4$, ensures
that, so far as the calculations concerning the Selmer group are
concerned, the prime $p=2$ behaves like a good prime.~\footnote{Top gives
condition (1a) as $a\equiv 3\pmod 4\text{ and }b\equiv 1\pmod 2\,
(\text{or }a\equiv 1\pmod 2\text{ and } b\equiv 2\pmod 4)$, but an
application of Tate's algorithm 
shows that the correct condition is given above.} Condition (1b), $a\equiv
2\pmod 3$, ensures that the prime $p=3$ also behaves like a good prime.  
\end{Rem}

The result of Top's theorem is an inequality for two reasons.  The first
reason is the unfortunate presence of the Tate-Shafarevich group.
This problem will not be dealt with here.  
The second reason that Top's theorem gives an inequality is that he merely
bounds the dimension of the Selmer groups by formulating the dimensions of
the groups into which he injects $S^\psi$ and $S^{\psi'}$.  Here we
make explicit the added restrictions necessary to determine the images of
the Selmer groups in the groups considered by Top. 

Because the primes 2 and 3 create special problems, the general result
will actually be stated for a generalized Selmer group which avoids these
primes.  The generalized Selmer group will agree with the standard Selmer
group when certain congruence conditions are placed on the parameter $a$,
and these conditions will be stated explicitly. 

\section{An Exact Sequence}\label{S:exact}

Let $\E$ be the elliptic curve considered above.  
Such an elliptic curve has a rational subgroup $T$ of order 3
generated by the point $(0,b\sqrt a)$.  Dividing by this subgroup yields
another elliptic curve $\Ep$, which has a rational subgroup $T'$ of order 3
generated by $(0,3(4a+27b)\sqrt{-3a})$.  Let $\psi:E\to E'$ be the quotient
map. Its dual, $\psi':E'\to E$, results from division by $T'$. 

Consider the fields $K=\Bbb Q(\sqrt{-3a})$ 
and $K'=\Bbb Q(\sqrt a)$. 
Following the notational conventions of Top, let $S$ denote the set of all
split primes of $K/\Bbb Q$. Write $N$ for the
map induced by the 
norm map, $N:K^*/{K^*}^3\to \Bbb Q^*/{\Bbb Q^*}^3$.  For a set of
primes~$\primes\in\Bbb Z$
which 
all split in $K$, write
\begin{equation}\label{E:hdef}
H(p_1,\dots,p_t)\buildrel\rm def\over=\{x\in\ker N:3|v_{\frak
p}(x)\,\forall\,\,\frak p\in S \text{ not lying over any }\primes\}.
\end{equation}
Likewise let $S'$, $N'$, and $H'$ be the corresponding objects for the
field
$K'$.

Top showed that $S^\psi$ injects into $\hset$ 
and $S^{\psi'}$ injects into $\hsetpr$, where
\begin{equation}\label{E:pdef}
P=\{\text{integer primes } p\ge 5 \text{ such that } 
p|(4a+27b) \text{ and }
(\frac{-3a}p)=1\}
\end{equation}
and
\begin{equation}\label{E:prdef}
P'=\{\text{integer primes } p\ge 5 \text{ such that }
p|b \text{ and }
(\frac{a}p)=1\}.
\end{equation}
We use a well-known exact sequence to extend these
facts. The following
lemma will be our key tool.  Its proof is a standard but 
tedious diagram chase.

\begin{Lem}\label{L:kercok}  
For any sequence of homomorphisms of abelian groups $A\buildrel f\over
\to B\buildrel g\over\to C$, we obtain the following exact sequence:
$$0\to \ker(f)\to\ker (g\circ f)\to\ker(g)\to\cok(f)\to\cok(g\circ
f)\to\cok(g)\to 0.$$
\end{Lem}

For every prime  $p\in\Bbb Q$ including $p=\infty$, we have the following
commutative diagram.
$$
\CD
0 @>>> E'(\Bbb Q)/\psi E(\Bbb Q) @>>> H^1(\Bbb Q,T) @>>> 
H^1(\Bbb Q,E)_\psi @>>> 0\\
@. @VVV @VVV @VVV @.\\
0 @>>> E'({\Bbb Q _p})/\psi E({\Bbb Q _p}) @>>> H^1({\Bbb Q _p}
,T) @>>> H^1({\Bbb Q _p},E)_\psi @>>> 0\\
@. @VVV @VVV @VVV @.\\
0 @>>> E'(\Bbb Q _p^{un})/\psi E(\Bbb Q _p^{un}) @>>>
H^1(\Bbb Q_p^{un},T) @>>> H^1(\Bbb Q_p^{un},E)_\psi @>>> 0
\endCD
$$

The commutative diagram is usually considered where the second and third
rows are taken as the direct sum over all primes $p$.  When we write
$\bigoplus_p$ with no conditions on $p$, we mean that the sum is to be
taken over all primes.  We also consider the diagram with second and
third rows taken as the direct sum over all primes $p\ne 2,3$.  This
allows us to avoid special problems which occur for these primes.

The portion of the diagram to which we apply Lemma~\ref{L:kercok} is 
\begin{equation}
H^1(\Bbb Q,T)\to\bigoplus_p H^1({\Bbb Q
_p},E)_\psi\to\bigoplus_p
H^1(\Bbb Q_p^{un},E)_\psi.
\end{equation}
Recall that $S^\psi$ is defined
to be the kernel of the first map.  Also note that $H^1(\Bbb R,E)$ is of
order dividing 2.  Since the order of $H^1(\Bbb R,E)_\psi$ divides both 2
and 3, the group must be trivial.  Therefore, when considering the above
diagram we may restrict our attention to the finite primes. 

We use the following notational convention.

\begin{Def}\label{D:gensel}
The {\it generalized Selmer group} is the kernel of
$$
H^1(\Bbb Q,T)\to\bigoplus_{p\ne 2,3} H^1({\Bbb Q _p},E)_\psi.
$$
It will be denoted $S_g^\psi$.
\end{Def}

Since we are interested in the dimension of the Selmer group, the
following easy lemma will
also be useful.
We will apply this lemma to the exact sequence given in Lemma~\ref{L:kercok}

\begin{Lem}\label{L:vectdim}
If $0\to V_1\to V_2\to\dots\to V_n\to 0$ is an exact
sequence of finite 
dimensional vector spaces, then
$$
\dim V_1-\dim V_2+\dim V_3-\dots+(-1)^{n-1}
\dim V_n=0.
$$ 
\end{Lem}



\section{The Middle Kernel of the Sequence}\label{S:middle}

We first analyze the kernel of $H^1(\Bbb Q,T)\to\bigoplus_{p} H^1(\Bbb
Q_p^{un},E)_\psi$.

\begin{Pptn}\label{P:ker1}
The kernel of $H^1(\Bbb Q,T)\to\bigoplus_{p\ne 2,3} H^1(\Bbb
Q_p^{un},E)_\psi$ is $\hset$, where $P=\{\text{primes
} p\ge 5 \text{ such that }p|(4a+27b) \text{ and }(\frac{-3a}p)=1\}$. 
\end{Pptn}

\emph{Proof.}
In this proof $p\ne 2,3$ is a prime.
In~\cite{jT93}, Top demonstrated
that $H^1(\Bbb Q,T)$ is isomorphic to $\ker N$, where
$K=\Bbb Q(\sqrt{-3a})$ and $N:K^*/{K^*}^3\to\Bbb Q^*/{\Bbb Q^*}^3$.
The isomorphism is well-determined up to sign, depending on the choice of
a generator for $T(\bar\Bbb Q)$.
Now $x\in\ker N$ implies that $3|v_p(N(x))$ for all $p$.  If we let $\fp$
denote a prime lying over $p$, then
\begin{equation}\label{E:normdiv}
v_p(N(x))=
\begin{cases}
v_\fp(x),&\text{if $p$ ramifies in $K$}\\
2v_\fp(x),&\text{if $p$ inert in $K$},
\end{cases}
\end{equation}
implies that $3|v_\fp(x)$ for nonsplit $p$.

First we consider the primes that ramify in $K$.  Let $\tilde{K}_\fp$
denote
the maximal unramified extension of $K_\fp$.  We have
\begin{equation}
H^1(\Bbb Q_p^{un},T)\cong\ker(Nm:\tilde{K}_\fp^*/\tilde{K}_\fp^{*3}\to
{\Bbb
Q_p^{un}}^*/{{\Bbb Q_p^{un}}^*}^3),
\end{equation}
where $\tilde{K}_\fp^*/\tilde{K}_\fp^{*3}\cong\Bbb Z/3\Bbb Z$ and ${\Bbb
Q_p^{un}}^*/{{\Bbb Q_p^{un}}^*}^3\cong\Bbb Z/3\Bbb Z$.
Since $v_p(N(x))=v_\fp(x)$, the norm map commutes with the identity map
on $\Bbb Z/3\Bbb Z$.  Therefore, $H^1(\Bbb Q_p^{un},T)=0$, and so
$H^1(\Bbb Q,T)\to H^1(\Bbb Q_p^{un},T)$ is the zero map.

If $p$ is not ramified, then $T\cong \mu_3$ over $\Bbb Q_p^{un}$.
Therefore, $\tilde{K}_\fp=\Bbb Q_p^{un}$, and so $H^1(\Bbb
Q_p^{un},T)\cong
{\Bbb Q_p^{un}}^*/{{\Bbb Q_p^{un}}^*}^3\cong \Bbb Z/3\Bbb Z$.  

If we identify $H^1(\Bbb
Q,T)$ with $\ker N$, it follows that for $p$ unramified, $x\in
H^1(\Bbb Q,T)$ has $3|v_\fp(x)$ if and only if it maps to the
identity in $H^1(\Bbb Q_p^{un},T)$.  So if $p$ is inert, then
$H^1(\Bbb Q,T)\to H^1(\Bbb Q_p^{un},T)$ is the zero map
by~\eqref{E:normdiv}.

Finally we treat the split primes.  By definition the elements of $\hset$
map to the identity in $H^1(\Bbb Q_p^{un},T)$ for all split primes
$p\notin P$.  Following the argument analogous to case 2 of page 310 of
Top~\cite{jT93}, we can show that $E'(\Bbb Q_p^{un})/\psi E(\Bbb
Q_p^{un})\cong \Bbb Z/3\Bbb Z$ for $p\in P$, in virtue of the following
commutative diagram with exact rows and columns.

$$
\CD
0 @>>> 0 @>>> 0\\
@VVV @VVV @VVV\\
E_0(\Bbb Q_p^{un}) @>>> E'_0(\Bbb Q_p^{un}) @>>>
 E'_0(\Bbb Q_p^{un})/\psi E_0(\Bbb Q_p^{un})\\
@VVV @VVV @VVV\\
E(\Bbb Q_p^{un}) @>>> E'(\Bbb Q_p^{un}) @>>>
 E'(\Bbb Q_p^{un})/\psi E(\Bbb Q_p^{un})\\
@VVV @VVV @VVV\\
\Bbb Z/n\Bbb Z @>>> \Bbb Z/3n\Bbb Z @>>> \Bbb Z/3\Bbb Z\\  
@VVV @VVV @VVV\\
0 @>>> 0 @>>> 0
\endCD
$$

The second-to-last row comes from Tate's algorithm, which shows that if
$p^n||(4a+27b)$, then the reduction for $E$ is type $\text{I}_n$, while
the reduction for $E'$ is type $\text{I}_{3n}$.  Since $E'_0(\Bbb
Q_p^{un})/\psi E_0(\Bbb Q_p^{un})=0$, we have that $E'(\Bbb
Q_p^{un})/\psi E(\Bbb Q_p^{un})\cong \Bbb Z/3\Bbb Z$. 

Thus, the exactness of 
\begin{equation}\label{E:bottomrow}
0\to E'(\Bbb Q_p^{un})/\psi E(\Bbb Q_p^{un})\to
H^1(\Bbb Q_p^{un},T)\to H^1(\Bbb Q_p^{un},E)_\psi\to 0
\end{equation}
implies that
$H^1(\Bbb Q_p^{un},E)_\psi=0$ for $p\in P$.  Therefore, $\hset$ is
contained in the kernel of $H^1(\Bbb Q,T)\to\bigoplus_{p\ne 2,3} H^1(\Bbb
Q_p^{un},E)_\psi$.

For the split primes at which $E, E'$ have good reduction, $E'(\Bbb
Q _p^{un})/\psi E(\Bbb Q
_p^{un})$ is the trivial group.  The remaining split primes are those
that divide $b$.  Using the argument analogous to case 2 of page 311 of
Top~\cite{jT93}, we can show that $E'(\Bbb
Q _p^{un})/\psi E(\Bbb Q
_p^{un})=0$ for these primes as well.  As before, we have a commutative
diagram.

$$
\CD
0 @>>> 0 @>>> 0\\
@VVV @VVV @VVV\\
E_0(\Bbb Q_p^{un}) @>>> E'_0(\Bbb Q_p^{un}) @>>>
 E'_0(\Bbb Q_p^{un})/\psi E_0(\Bbb Q_p^{un})\\
@VVV @VVV @VVV\\
E(\Bbb Q_p^{un}) @>>> E'(\Bbb Q_p^{un}) @>>>
 E'(\Bbb Q_p^{un})/\psi E(\Bbb Q_p^{un})\\
@VVV @VVV @VVV\\
\Bbb Z/3n\Bbb Z @>>> \Bbb Z/n\Bbb Z @>>> 0\\  
@VVV @VVV @VVV\\
0 @>>> 0 @>>> 0
\endCD
$$

The second-to-last row comes from Tate's algorithm, which shows that if
$p^n||b$, then the reduction for $E$ is type $\text{I}_{3n}$, while
the reduction for $E'$ is type $\text{I}_n$
Thus we have an isomorphism $ E'_0(\Bbb Q_p^{un})/\psi E_0(\Bbb
Q_p^{un})\cong E'(\Bbb Q_p^{un})/\psi E(\Bbb Q_p^{un})$.  Since the
former quotient is zero, so is the latter.
  
Thus, for the split primes $p\notin P$, the exactness of
~\eqref{E:bottomrow} implies that an element of $H^1(\Bbb Q,T)$ will map
to the identity in $H^1(\Bbb Q_p^{un},E)_\psi$ if and only if it maps to
the identity in $H^1(\Bbb Q_p^{un},T)$. Therefore, the kernel of
$H^1(\Bbb Q,T)\to\bigoplus_{p\ne 2,3} H^1(\Bbb Q_p^{un},E)_\psi$ is
precisely $\hset$.\qed

If we wish a result for the standard Selmer group, we must put conditions
on $a$ to make things work correctly at $p=2,3$. 

\begin{Pptn}\label{P:ker1g2}
The kernel of $H^1(\Bbb
Q,T)\to\bigoplus_{p\ne 3} H^1(\Bbb Q_p^{un},E)_\psi$ is $\hset$,
when $a\equiv 2\text{ or }3\pmod 4$.
\end{Pptn}

\emph{Proof.} If $a\equiv 2\text{ or }3\pmod 4$, then 2 ramifies in $K$,
so by the argument in the previous proof, $H^1(\Bbb Q,T)\to H^1(\Bbb
Q_p^{un},T)$ is the zero map.\qed

\begin{Def}
Denote by $H(P)_r$ the  subset
of $\hset$ represented by elements $x$ such that
$x\in \tilde{K}_\fp^{*3}$ for a choice of $\fp|3$.  
\end{Def}

\begin{Pptn}\label{P:ker1g3} 
When $a\equiv 2 \text{ or } 3\pmod {4}$ and $3\nmid a$, the kernel of the
map $H^1(\Bbb Q,T)\to\bigoplus_p H^1(\Bbb Q_p^{un},E)_\psi$ is 
$H(P)_r$.  
\end{Pptn}

\emph{Proof.}  The condition that $3\nmid a$ implies that 3
ramifies in $K$.  As in the previous theorem, we have 
\begin{equation}
H^1(\Bbb Q_3^{un},T)\cong\ker(Nm:\tilde{K}_\fp^*/\tilde{K}_\fp^{*3}\to
{\Bbb
Q_3^{un}}^*/{{\Bbb Q_3^{un}}^*}^3).
\end{equation}
Using Tate's
algorithm~\cite{jS94}, we find that $E'$ has reduction of type
$\text{I}_{v_3(b)}$, and $E$ has reduction of type $\text{I}_{3v_3(b)}$.
Therefore,
$E'(\Bbb Q_p^{un})/\psi E(\Bbb Q_p^{un})=0$, and so the kernel of
$H^1(\Bbb Q,T)\to H^1(\Bbb Q_3^{un},E)_\psi$ is the same as the kernel of
$H^1(\Bbb Q,T)\to H^1(\Bbb Q_3^{un},T)$.  Using the 
explicit description
of $H^1(\Bbb Q_3^{un},T)$ and the previous
propositions, we find that the
kernel of  $H^1(\Bbb Q,T)\to\bigoplus_p H^1(\Bbb Q_p^{un},E)_\psi$ is
$H(P)_r$.
\qed

Combining Proposition~\ref{P:ker1g3} with Lemma~\ref{L:kercok}, we
obtain Top's result
that $S^\psi$ maps into $\hset$.  Since we are interested in the
dimension of $S^\psi$, we would like to know the dimension of $\hset$. 
The following proposition gives a formula for the dimension of this
vector space. 

Denote the elements of the set $P$ by $p_1,\cdots,p_t$.  By definition,
these primes split in the field $K=\Bbb Q(\sqrt{-3a})$.  Choose a set of
primes of $K$, $\{\fp_1,\cdots,\fp_t\}$, so that $\fp_i|p_i$.

\begin{Pptn}\label{P:hpdim}
Let $V=\{(i_1,\dots,i_t)\in\zzt\, \text{\rm such that }\,
\fp_1^{i_1}\dots\fp_t^{i_t}\in
Cl(K)^3\}$, then
$$\dim_{\mathbf F_3}\hset=r_3(-3a)+\dim_{\mathbf F_3}V+\dim_{\mathbf
F_3}U/U^3,$$ where $U$ denotes the units in the ring of integers of
$K$.
\end{Pptn}

\emph{Proof.}
Define the map $\Theta$ by
\begin{align}
\Theta:\hset&\longrightarrow\zzt,\\
x\qquad&\longmapsto\bigoplus_{i=1}^t\, v_{\fp_i}(x)\pmod 3,
\end{align}
where $\fp_i\subseteq\cal O_K$ is one of the primes lying over $p_i$.
Suppose 
$x\in\ker\Theta$.  
It is easy to see that the ideal of $K$ generated by $x$ must
be a cube.
Therefore the
dimension of $\ker\Theta$ is $r_3(-3a)+\dim_{\mathbf F_3}U/U^3$. 

The following lemma gives a formulation for $\image\Theta$. 
Lemma~\ref{L:vectdim} then gives the desired equality.\qed

\begin{Lem}\label{L:imtheta}
The image of $\Theta$ is the vector space $V$.
\end{Lem}

\emph{Proof.}
If $(i_1,\dots,i_t)\in\image\Theta,$ then there is an $x\in\hset$ with the 
ideal factorization
$(x)=\idfac ik\frak b^3$ for some 
$\frak b\subseteq\cal O_K$ where $k_j=0,\,1,\text{ or }2$ according as 
$i_j=0,\,2, \text{ or }1$.  Since the $p_i$ split, we have an equality of ideals 
$(p_1^{i_1}\dots p_t^{i_t})=\idfac ii$.  
Because $i_j+k_j\equiv 0\pmod 3$,  
multiplying the ideal equality by the factorization for $(x)$ gives
$(p_i^{i_1}\dots p_t^{i_t}x)=\idfacsh {2i}\frak a^3$ for a different 
ideal $\frak a$.
Squaring both sides gives $(y)=
\idfacsh i\frak c^3$ for still another ideal $\frak c$, where $y=p_1^{2i_1}\dots
p_t^{2i_t}x^2$.  This implies that
$\idfacsh i\sim(\frak c^{-1})^3\text{ in }Cl(K)$, or in other words 
$(i_1,\dots,i_t)\in V$.

On the other hand, if $(i_1,\dots,i_t)\in V$, then  there is a $y\in K$
with
the property that $(y)=\idfacsh i\frak b^3$.  Multiplying by the ideal
$(p_1^{i_1}\dots p_t^{i_t})=\idfac ii$ and squaring yields $(x)\idfac
i{2i}\frak a^3$ for some ideal $\frak a$, where $x=p_1^{2i_1}\dots
p_t^{2i_t}y^2$.  The element $x$ satisfies the conditions to be in
$\hset$, and $\Theta(x)=(i_1,\dots,i_t)$.\qed 

\begin{Rem}
Since the dimension of $V$ is clearly at most $t$,
Proposition~\ref{P:hpdim} implies
Lemma 3 of~\cite{jT93}, in which 
Top
gives an inequality for the dimension of $\hset$. 
\end{Rem}

\begin{Rem} {\it A priori\/} to determine if the vector $(i_1,\dots,i_t)$
is in the image of $\Theta$, one must search for any $x\in K^*$ with the
proper ideal factorization.  The factorization specifies only what must
happen at the primes lying over $\primes$.  Thus, the vector will be in
the image if any one of the infinitely many allowed prime factorization
gives a principal ideal.  Lemma~\ref{L:imtheta} replaces this infinite
search with the test of a single ideal $\idfacsh i$.  \end{Rem}


\section{The Last Kernel of the Sequence}\label{S:last}

Let $\E$ be the usual elliptic curve without congruence conditions on $a$
or $b$.
In the
fundamental exact sequence with which we work, the third group is the
kernel of the map 
\begin{equation}
\bigoplus_p H^1(\Bbb Q_p,E)_\psi\to\bigoplus_p
H^1(\Bbb
Q_p^{un},E)_\psi.
\end{equation}
We analyze
this kernel by first considering the kernel of a single factor of this 
map. 

Let $G=
\gal(\overline{\Bbb Q_p}/\Bbb Q_p)$, and $I=
\gal(\overline{\Bbb Q_p}/\Bbb Q_p^{un})$.
By the inflation-restriction sequence, 
\begin{equation}
0\to H^1(G/I 
,E(\Bbb Q_p^{un}))\to H^1(\Bbb Q_p,E)\to H^1(\Bbb 
Q_p^{un},E)^{G/I}\to \dots.
\end{equation}
Furthermore, Milne shows in Proposition I.3.8 of~\cite{jsM86} that 
\begin{equation}
H^1(G/I 
,E(\Bbb Q_p^{un}))\cong H^1(\zh,\pio),
\end{equation}
where $\pio=E(\Bbb Q_p^{un})/E_0(\Bbb Q_p^{un})$, and 
$\gal(\Bbb Q_p^{un}/\Bbb Q_p)$ is canonically
identified with $\zh$.  
Taking $\psi$-kernels of all the cohomology groups and putting together the
factors for each $p$, we find that $\hset$ maps into 
$\bigoplus_p H^1(\hat{\Bbb Z},\pi_0(\overline E))_\psi$ in our
exact sequence.

\begin{Pptn}\label{P:ker2} 
The kernel of $\,\bigoplus_{p\ne 2,3} H^1(\Bbb
Q_p,E)_\psi\to\bigoplus_{p\ne 2,3} H^1(\Bbb Q_p^{un},E)_\psi$ is
$\bigoplus_{p\in P'}H^1(\zh,\pio)_\psi\cong(\Bbb Z/3\Bbb Z)^{\#P'}$ where
$P'$ is the set of primes $\{p\ge 5\text{\rm { }such that } p|b\text{\rm {
}and }(\frac ap)=1\}$.  Moreover, if $a\equiv 2 \text{ or }
11\pmod {12}$, then this is also the kernel of $\,\bigoplus_p H^1(\Bbb
Q_p,E)_\psi\to\bigoplus_p H^1(\Bbb Q_p^{un},E)_\psi$.  
\end{Pptn}

\emph{Proof.} Top~\cite{jT93} gives the possible reduction types of $E$
and $E'$ for $p\ge 5$ to be $\text{I}_0,\,\text{II},\,\text{I}_n^*, \text{
and }\text{I}_n$.  Using Tate's algorithm~\cite{jS94} we see that with
$a\equiv 2,3\pmod 4$, the possible reduction types
at $p=2$ are $\text{II}\text{ and I}_n^*$. 
We can also calculate that $E$ has
reduction of type $\text{I}_{3v_3(b)}$ at $p=3$.

First we deal with the cases $p\ne 3$.  If the reduction is of type
$\text{I}_n^*$,
the Kodaira-N\'eron classification gives that 
\begin{equation}
\pio\cong
\begin{cases}
 \Bbb Z/2\Bbb Z\times\Bbb Z/2\Bbb Z,&\text{for $n$ even}\\
\Bbb Z/4\Bbb Z,&\text{for $n$ odd}.
\end{cases}
\end{equation}
Since $\psi\circ\psi'$ is multiplication by 3, 
$\ker\psi\subseteq\ker{[3]}$.  Because $\pio$ has no elements of order 3,
its 
$\psi$-kernel is trivial.  Therefore, $\psi:\pio\to \pi_0(\overline
E')$ defines an isomorphism,
which implies that $H^1(\zh,\pio)_\psi$ is also 0.
 
If the reduction is of type $\text{I}_0 \text{ or II}$, then
$H^1(\zh,\pio)_\psi=0$, since $\pio=0$.

Finally, if the reduction is of type $\text{I}_n$, then $\pio$ is cyclic.  
The action of $\zh$
on $\pio$ must preserve the structure, so $\zh$ must
act as $\pm 1$, as in Theorem 3.7 of~\cite{eS96}.
Consider the short exact sequence
\begin{equation}
0\to U\to \zh\to \{\pm 1\}\to 0,
\end{equation}
where $U$ is the kernel of the map sending an element of $\zh$ to its action.
Then $U$ acts trivially on $\pio$.  Using the inflation-restriction
sequence, we obtain
\begin{multline}
0\to H^1(\{\pm 1\},\pio)\to H^1(\zh,\pio)\to \\
H^1(U,\pio)^{\{\pm 1\}}\to H^2(\{\pm
1\},\pio)\to\dots.
\end{multline}
By counting, we find that both $H^1(\{\pm 1\},\pio)_\psi$ and $H^2(\{\pm
1\},\pio)_\psi$ are 0.  Thus, by left exactness 
$H^1(\zh,\pio)_\psi= H^1(U,\pio)^{\{\pm 1\}}_\psi$.

Since $U$ acts trivially and has a single topological generator, we find
that 
\begin{equation}
H^1(U,\pio)^{\{\pm
1\}}\cong\operatorname{Hom}{(U,\pio)^{\{\pm
1\}}}\cong\pio^{\{\pm 1\}}.
\end{equation}
Looking at $\psi$-kernels, we obtain
\begin{equation}
H^1(\zh,\pio)_\psi\cong\pio_\psi^{\{\pm 1\}}.
\end{equation}

Top shows that, in the case of multiplicative reduction, $p|b$ implies
that $\pio_\psi\cong\Bbb Z/3\Bbb Z$ and $\pio_\psi=0$ otherwise.  Because
$\zh$ acts on $\pio_\psi$ as it acts on $\pio$, $H^1(\zh,\pio)_\psi$ will
give a copy of $\Bbb Z/3\Bbb Z$ when $\zh$ acts as $+1$, and the trivial
group otherwise. 

By IV.9.6 of~\cite{jS94}, $\zh$ acts on $\pio$ as $+1$ if $E$ has split
multiplicative reduction and as $-1$ if $E$ has non-split multiplicative
reduction.  Since $p|b$, the reduced equation is $\tilde E:\tilde
y^2=\tilde x^3+a\tilde x^2$.  Thus, the reduction is split precisely when
$(\frac ap)=1$.

Finally, we consider the case $p=3$.  The condition that $a\equiv 2\pmod
3$ implies that the splitting field of $T^2-a$ over ${\bold F}_3$
strictly contains ${\bold F}_3$.  Therefore, Tate's algorithm implies
that $\zh$ acts as $-1$, and so $\pio^{\{\pm 1\}}$ has order 1 or 2,
which means that its $\psi$-kernel is trivial.\qed

\section{A Formula for $S^\psi$}\label{S:formula}

By previous calculations, the generalized Selmer group fits into
the exact sequence
\begin{equation}
0\to S_g^\psi\to\hset\to\bigoplus_{p\in P'} (\Bbb Z/3\Bbb Z)\to\dots,
\end{equation}
where the sets $P$ and $P'$ are defined in Section 2.2.
Therefore, $S_g^\psi$ is isomorphic to the kernel of $\hset\to
\bigoplus_{p\in P'}
(\Bbb
Z/3\Bbb Z)$.  Moreover, if $a\equiv 2 \text{ or } 11\pmod {12}$, then the
standard Selmer
group
fits into the exact sequence
\begin{equation}
0\to S^\psi\to H(P)_r\to\bigoplus_{p\in P'} (\Bbb Z/3\Bbb Z)\to\dots,
\end{equation}
where $H(P)_r$ denotes the subset of $\hset$ given in
Proposition~\ref{P:ker1g3}.
Thus, $S^\psi\cong\ker(H(P)_r\to\bigoplus_{p\in P'} (\Bbb Z/3\Bbb Z))$.
To determine these kernels we will analyze the map
$\hset\to\Bbb
Z/3\Bbb Z$ for each prime $p\in P'$.  We  allow $P'$ to
denote the aforementioned set of integral primes as well as the set of
primes lying over these in any field extension.

\begin{Lem}\label{L:pprimeisom}
For $p\in P'$, $H^1(\Bbb Q_p,T)\cong H^1(\Bbb Q_p,E)_\psi$.
\end{Lem}

\emph{Proof.}
As in the proof of~\ref{P:ker1}, the application of Tate's algorithm for
primes $p|b$ shows that $E'(\Bbb Q_p)/\psi E(\Bbb Q_p)$ is trivial.  By
the exactness of
\begin{equation}
0\to E'({\Bbb Q _p})/\psi E({\Bbb Q _p})\to H^1({\Bbb Q _p},T)
\to H^1({\Bbb Q _p},E)_\psi\to 0,
\end{equation}
the result follows.\qed

If $x$ is in the kernel of $\hset\to \bigoplus_{p\in P'}(\Bbb Z/3\Bbb Z)$,
then it must map to the identity in $H^1(\Bbb Q_p,E)_\psi$ for all $p\in
P'$.  By the previous lemma this implies that $x$ maps to the identity in
$H^1(\Bbb Q_p,T)$ for these primes.  Therefore, $x$ is in the kernel if it
maps to a cube in $K_\fp^*$ for all $\fp\in P'$.  By the usual
decomposition $K_\fp^*\cong\fp^{\Bbb Z}\times U_{K_\fp}$, we see that
${K_\fp^*}^3\cong\fp^{3\Bbb Z}\times {U_{K_\fp}}^3$.  By the definition of
$\hset$, $v_\fp(x)\equiv 0\pmod 3$ for all $\fp\in P'$,
since $P$ and $P'$ are disjoint. Let 
\begin{equation}
u=\frac x{\fp^{v_\fp(x)}}.
\end{equation}
 Then,
$u\equiv x\pmod {{K^*}^3}$ and $v_\fp(u)=0$.  For the coset represented by
$x$ to be in the Selmer group, it must be that $u\in U_{K_\fp}^3$. We have
proved the following two theorems. 

\begin{Thm}\label{T:Sg}
The generalized Selmer group $S^\psi_g$ is the subgroup of $\hset$
represented by elements with
$v_\fp(u_\fp)=0$ and $u_\fp\in {U_{K_\fp}}^3$ for all
$\fp\in P'$.
\end{Thm}

\begin{Thm}\label{T:S}
If $a\equiv 2\text{ or }11\pmod{12}$, then $S^\psi$  is the subgroup of
$\hset_r$
represented by elements with
$v_\fp(u_\fp)=0$ and $u_\fp\in {U_{K_\fp}}^3$ for all
$\fp\in P'$.
\end{Thm}

Let us name the maps
\begin{equation}
\Phi:\hset\to\bigoplus_{p\in P'}\Bbb Z/3\Bbb Z
\end{equation}
and
\begin{equation}
\hat\Phi:\hset_r\to\bigoplus_{p\in P'}\Bbb Z/3\Bbb Z,
\end{equation}
which send a coset to the direct sum of its unit representatives
$u_\fp\pmod{U_{K_\fp}^3}$.

\begin{Cor}\label{c:dimS}
The dimension of the generalized Selmer group is
$$
\dim_{\mathbf F_3}S_g^\psi=r_3(-3a)+\dim_{\mathbf F_3}V+\dim_{\mathbf
F_3}U_K/U_K^3-\dim_{\mathbf F_3}\image\Phi,$$
where $V$ is the vector space defined in Proposition~\ref{P:hpdim}.  

Moreover, if $a\equiv
2\text{ or }11\pmod{12}$, then the dimension of the Selmer group is
$$\dim_{\mathbf F_3}S^\psi=r_3(-3a)+\dim_{\mathbf F_3}V+\dim_{\mathbf
F_3}U_K/U_K^3-\dim_{\mathbf F_3}\image\hat\Phi-\nu,$$
where $\nu$ is the codimension of $\hset_r$ in $\hset$.
\end{Cor}

\emph{Proof.}
Since $S_g^\psi=\ker\Phi$ and $S^\psi=\ker\hat\Phi$, this follows immediately
from Lemma~\ref{L:vectdim} and Proposition~\ref{P:hpdim}.\qed

\begin{Rem}
We note that similar results for the Selmer group of $\psi'$ follow, with
the sets of primes $P$ and $P'$ trading roles.  The only changes are in
the proof of the
analogue to Proposition~\ref{P:ker1g3}, where we note that 3 is inert in
$Q(\sqrt{a})$ when $a\equiv 2\pmod 3$, and in the proof of the analogue to
Proposition~\ref{P:ker2}, where we note that the possible reduction types
for $E'$
are the same as for $E$ when $a\equiv 2\pmod 3$.
\end{Rem}


\section{Relating the Selmer Groups of Dual Isogenies}\label{S:relate}

In this section, we directly relate the dimension of the Selmer group
$S^\psi$ to the dimension of $S^{\psi'}$ using a formula of
Cassels~\cite{jC65}.  
The formula is
given by
\begin{equation}\label{E:Cassels}
\frac{\#S^\psi}{\#S^{\psi'}}=\frac{\#T(\Bbb Q)}{\#T'(\Bbb
Q)}\frac{\alpha'}{\alpha}\prod_p\frac{c_p'}{c_p},
\end{equation}
where $\alpha=\int_{E(\Bbb R)}\omega$ (resp. $\alpha'=\int_{E'(\Bbb
R)}\omega'$) for $\omega$ the canonical N\'eron differential of $E$ (resp.
$\omega'$ the canonical N\'eron  differential of $E'$), $c_p$ (resp.
$c_p'$) is the number of connected components of the special fiber of the
N\'eron model of $E$ (resp. $E'$), and $T(\Bbb Q)$ (resp. $T'(\Bbb Q)$) is
the group of points of the kernel of $\psi$ (resp. $\psi'$) which are defined
over $\Bbb Q$.  Since $T(\Bbb Q)$ and $T'(\Bbb Q)$ are trivial, to relate
the dimensions of the Selmer groups we only need to 
find the reduction types and calculate the integrals of the canonical
differentials.  Throughout this section we consider only $a\equiv 2\text{
or }11\pmod{12}$.

First we calculate $c_3$ and $c_3'$.  By assumption 3 does not divide
$a$, and therefore does not divide $(4a+27b)$.   
Since 
\begin{equation}
\Delta_E=-2^4a^2b^3(4a+27b),
\end{equation}
we see that $v_3(\Delta_E)=3v_3(b)$.  Using Tate's algorithm for $p=3$ we
find that, since $3\nmid b_2=4a$, the reduction is of type
$\text{I}_{3v_3(b)}$, and so
$c_3=3v_3(b)$.
On the other hand, 
\begin{equation}
\Delta_{E'}=-2^43^{12}a^2b(4a+27b)^3.
\end{equation}
Therefore the equation for $E'$ is not minimal at $p=3$.  It becomes
minimal at 3 after the
change of
variables 
\begin{equation}\label{E:changexy}
x=3^2x',\quad y=3^3y'.
\end{equation}
The new discriminant is $2^4a^2b(4a+27b)^3$.  An application of
Tate's algorithm 
then yields that the reduction is of
type $\text{I}_{v_3(b)}$, and so $c_3'=v_3(b)$.

Next we calculate $c_p$ and $c_p'$ for $p\ne 3$.  For $p=2$ Tate's
Algorithm shows that $c_2=c_2'$.  For $p\ge 5$, Top~\cite{jT93} showed
that the reduction is of type
$\text{I}_0,
\text{II}, \text{ or I}_n^*$, and hence $c_p=c_p'$, unless either 
\begin{equation}
p\nmid a \text{ and }p|b
\end{equation}
or
\begin{equation}
p\nmid a \text{ and }p|(4a+27b).
\end{equation}
In the former case, $c_p=3v_p(b)$
and
$c_p'=v_p(b)$.  In the latter case, $c_p=v_p(4a+27b)$
and
$c_p'=3v_p(4a+27b)$.

We now calculate the ratio $\alpha'/\alpha$.  The explicit formula for
the map $\psi$ is given in~\cite{jT93} by
\begin{equation}
\psi(x,y,z)=(9(2xy^2z+2ab^2xz^3-x^4-\frac{2}{3}ax^3z),27y(4abxz^2-8ab^2z^3+x^3)
,x^3z).
\end{equation}
Therefore, using the quotient rule, we obtain
$$\psi^*\omega''=\frac 13\omega,$$
where $\omega$ is the invariant differential for the given Weierstrass
equation for $E$ and  $\omega''$ is the invariant differential for the
given Weierstrass
equation for $E'$.
Since the Weierstrass equation for $E'$ is not minimal at $p=3$, we use the 
change of variables~\eqref{E:changexy}, and calculate that
\begin{equation}
\frac{1}{3}\omega'=\omega''.
\end{equation}
Therefore, we find that
\begin{equation}
\psi^*\omega'=\omega.
\end{equation}
This relation and a little calculus gives
\begin{equation}
\alpha=\int_{E(\Bbb R)}\omega=\int_{E(\Bbb R)}\psi^*\omega'
=\int_{\psi(E(\Bbb R))}\omega'.
\end{equation}

The real loci of the curves $E$ and $E'$ are connected when
$b$ and $4a+27b$ have the same sign, but the curves each have two
connected
components when $b$ and $4a+27b$ have opposite signs.  
It is easy to verify that the connected component of $E(\Bbb R)$
containing infinity is a one-fold (resp. three-fold) cover of the
connected component of $E'(\Bbb R)$ containing infinity, when $a<0$ (resp.
$a>0$).  Likewise, if the second connected components are present, then
the component not containing infinity of  $E(\Bbb R)$
is a one-fold (resp. three-fold) cover of the
connected component of $E'(\Bbb R)$ not containing infinity when $a<0$
(resp.
$a>0$).
Therefore, 
\begin{equation}
\int_{\psi(E(\Bbb R))}\omega'=
\begin{cases} \int_{E'(\Bbb R)}\omega',&\text{if
$a<0$},\\
3\int_{E'(\Bbb R)}\omega',&\text{if $a>0$}.
\end{cases}
\end{equation}

Combining these results yields
\begin{equation} 
\alpha=
\begin{cases} \alpha',&\text{if $a<0$}\\
3\alpha',&\text{if $a>0$}.
\end{cases}
\end{equation}

The relationship~\eqref{E:Cassels} then gives us the following
formula.

\begin{Pptn}  The dimensions over $\mathbf F_3$ of the Selmer groups
$S^\psi$ and $S^{\psi'}$ are subject to the relation
$$
\dim_{\mathbf F_3}S^\psi-\dim_{\mathbf F_3}S^{\psi'}=A-B-\delta-\epsilon,
$$
where $A=\#\{p\mid p\ge 5, p\nmid a \text{ and }p\mid (4a+27b)\}$, 
$B=\#\{p\mid p\ge
5, p\nmid a \text{ and }p\mid b\}$, 
$$
\delta=
\begin{cases}
 0, &\text{if $3\nmid b$}\\
1,&\text{if $3|b$}.
\end{cases}
$$
and 
$$
\epsilon=
\begin{cases}
 0, &\text{if $a<0$}\\
1,&\text{if $a>0$}.
\end{cases}
$$
\end{Pptn}

We can compare this to the results obtained by subtracting the
dimensions of the Selmer groups obtained in the previous analysis.  Using
the formulas for the exact dimensions of the Selmer groups yields
\begin{multline}
\dim_{\mathbf F_3}S^\psi-\dim_{\mathbf
F_3}S^{\psi'}=\\
r_3(-3a)-r_3(a)+\dim_{\mathbf F_3}V-\dim_{\mathbf
F_3}V'-\dim_{\mathbf F_3}\hat\Phi+\dim_{\mathbf
F_3}\hat\Phi'+\beta+\nu-\nu'.
\end{multline}
Here  
\begin{equation}
\beta=
\begin{cases} 1, &\text{if $a<0$}\\
-1,&\text{if $a>0$}.
\end{cases}
\end{equation}
Note that Scholz's theorem implies that 
\begin{equation}
r_3(-3a)-r_3(a)=
\begin{cases} 0\text{ or }-1,&\text{if
$a<0$}\\
0\text{ or }1,&\text{if $a>0$}.
\end{cases}
\end{equation}
Also recall that $V$ and $\Phi'$ are vector spaces of dimension at most
\begin{equation}
C=\#\{p\ge 5\mid p|(4a+27b)\text{ and }\left(\frac{-3a}p\right)=1\},
\end{equation}
and that $V'$ and $\Phi$ are vector spaces of dimension at most
\begin{equation}
D=\#\{p\ge 5\mid p|b\text{ and }\left(\frac ap\right)=1\}.
\end{equation}
We note that, heuristically speaking, $C$ is roughly $A/2$ and $D$ is
roughly $B/2$.  Therefore, the differences obtained via the two methods are
consistent.  The difference obtained by utilizing the theorem of
Cassels is more concrete and much easier to compute.  However, appealing
to this process only allows the calculation of the difference of the
dimensions of the Selmer
groups and not the sum.  Finding an upper bound for the rank of the
elliptic curves
requires knowing the sum of the dimensions of the Selmer groups.
Therefore, the lengthier analysis provides information not available
through the duality theorem.